\documentclass[a4paper]{amsart}

\usepackage{amssymb}
\usepackage{amsthm}
\usepackage{amsmath}
\usepackage{amscd}
\usepackage{graphicx}
\usepackage{dsfont}
\usepackage{xcolor}
\renewcommand\a{\mathfrak{a}}

\newcommand\op[1]{\mathop{\rm #1}\nolimits}

\newcommand\R{\mathbb{R}}

\newcommand\g{\mathfrak{g}}
\newcommand\h{\mathfrak{h}}
\newcommand\La{\Lambda}
\newcommand\m{\mathfrak{m}}
\newcommand\ot{\otimes}
\renewcommand\th{\theta}
\newcommand\vp{\varphi}

\theoremstyle{plain}
\newtheorem{theorema}{Theorem}
\newtheorem{theoremb}{Theorem}
\newtheorem{theoremc}{Theorem}

\newtheorem{example}[theoremc]{Example}

\begin{document}
\title[Erratum: non-degenerate almost complex structures in 6D]{Erratum to:
Almost complex structures in 6D with non-degenerate Nijenhuis tensors\\ and large symmetry groups}
\author{B. Kruglikov, H. Winther}
 \address{Department of Mathematics and Statistics, Faculty of Science and Technology,
 UiT the Arctic University of Norway, 
 Troms\o\ 90-37, Norway}
 \email{ boris.kruglikov@uit.no, \quad henrik.winther@uit.no.}
 \maketitle

 \begin{abstract}
We correct an error in the second part of Theorem 3 of our original paper.
 \end{abstract}

\maketitle

It was stated in \cite[Theorem 3]{KW1} that nondegenerate almost complex manifolds $(M^6,J)$ with a symmetry
algebra of dimension 9 are all (locally) homogeneous spaces $M=G/H$ and that they have a semi-simple stabilizer $H$.
While the first part is correct, the second claim fails as the following shows.


 \begin{example}
 Let $\a=\R b \oplus \R z \oplus V^4$ be the solvable Lie algebra with $\R z\oplus V^4$ a 5D Heisenberg sub-algebra and $b$ acting as a derivation on $\R z\oplus V^4$ with weight $1$ on $V^4$ and $2$ on $\R z$. Note that $\a$ is isomorphic to the radical of the algebra $A.3.1$ from \cite{AKW} (for a particular choice of parameters). The outer derivations of $\a$ are given by
 $$
 \op{out}(\a):=\op{der}(\a)/\a=\mathfrak{sp}(4,\R),
 $$
and the action of $\mathfrak{sp}(4,\R)$ is standard on $V^4$ and trivial on $\R z \oplus \R b$. There are several inequivalent embeddings of $\mathfrak{sl}(2,\R)$ into $\mathfrak{sp}(4,\R)$. The two most interesting to us are the embeddings $\mathfrak{su}(1,1)$, which preserves a complex structure on $V^4$, and $\mathfrak{sl}_2^{irr}$, which acts irreducibly on $V^4$. If we let $\g=\mathfrak{su}(1,1)\ltimes \a$, then $\g$ is the algebra A.3.1 from \cite{AKW}. However, we let $\g=\mathfrak{sl}_2^{irr}\ltimes \a$ instead. Then $\g$ is no longer isomorphic to A.3.1, and $\g/\mathfrak{sl}_2^{irr}$ does not admit an $\mathfrak{sl}_2^{irr}$-invariant almost complex structure; but there is a 3D solvable sub-algebra $\h=\mathfrak{l}_1\subset \g$, which has a one-dimensional intersection with $V^4\subset\a$, such that
$$
\mathfrak{m}=\g/\h
$$
admits an $\h$-invariant almost complex structure, and this structure is non-degenerate. There is a basis $z,b,x_1,\ldots,x_4,h,e,f$ of $\g$ such that the structure equations of $\g$ are:
 \begin{gather*}
[b, x_i] =  x_i,\ [b,z] = 2z,\
[h,x_1] = -3 x_1,\ [h,x_2] = -x_2,\  [h,x_3] = x_3,\ [h,x_4] = 3 x_4, \\
[f,x_2] = -3 x_1,\ [f,x_3]= -2 x_2,\ [f,x_4] = -x_3,\
[e,x_1] = x_2,\ [e,x_2] = 2 x_3,\ [e,x_3]= 3 x_4,\\
[x_1,x_4] = z,\ [x_2,x_3] = - 3 z,\ [h, f] = -2 f,\ [h,e] = 2 e,\ [f, e] = h.
 \end{gather*}
In this basis, we have $\mathfrak{l}_1=\langle x_1, h -b, f -z\rangle$.
Identifying $\m=\g/\h=\langle z,b,e,x_2,x_3,x_4\rangle$ the almost complex
structure is given by $Jz=x_2,Jb=-x_3,Je=x_4$ (the invariant structure $J$ is unique: the conjugation $J\mapsto-J$ is obtained
by the outer automorphism $x_i\mapsto-x_i$).

We also note that there exists a unique (up to scale) invariant almost pseudo-Hermitian metric $(g,J)$ on $\m$.
In the dual basis $z^*,b^*,x_1^*,\ldots,x_4^*$: $g=e^*z^*+x_2^*x_4^*-b^{*2}-x_3^{*2}$.
 \end{example}

This counter-example was found upon revisiting the result via a new technique. In \cite{KW1} we relied
on a Maple computation with polynomial ideals that express the Jacobi identity for reconstructed Lie algebra
structure of the symmetry algebra $\g$. Some solutions have been lost with this approach. In \cite{KW2} we
elaborated a different reconstruction technique, which effectively separates linear constraints from
genuine quadratic relations, and exploring it we obtained the missing case(s), thus complementing the classification of homogeneous spaces (with semi-simple isotropy) performed in \cite{AKW}.

Recall from Section 3 of \cite{KW1} that the isotropy algebra $\h$ can be one of the types: $\mathfrak{p}$,
$\mathfrak{r}$, $\mathfrak{l}_0$, $\mathfrak{l}_1$ or $\mathfrak{l}_2$. In the cases $\h$ is $\mathfrak{p}$,
$\mathfrak{l}_2$ or $\mathfrak{r}$ (when $s\not\in\mathfrak{r}$) we computed
$H^1(\h,\op{Hom}(\m,\h))=0$, so the homogeneous space is reductive and the further computations hold.
We should only revisit the remaining cases.

 \begin{theorema}
Let $(G/H,J)$ be a homogeneous almost complex 6D manifold with non-degenerate Nijenhuis tensor $N_J$.
If the stabilizer $H$ is not semi-simple, then its Lie algebra is $\h=\mathfrak{l}_1$, and $\g=\mathfrak{sl}_2\ltimes \a$, with $\mathfrak{sl}_2$ acting irreducibly on $V^4 \subset \mathfrak{a}$ as described in the example above, and $\mathfrak{l}_1$ is conjugate to the one given in the example.
 \end{theorema}

Before proving this theorem, we recall the main result of \cite{KW2}. 
For $\vp\in\h^*\ot\m^*\ot\h$, $h\in\h$ and $u_1,u_2\in\m$ define
$\delta\vp\in\h^*\ot\La^2\m^*\ot\m$, $Q\vp\in\h^*\ot\La^2\m^*\ot\h$ by
 \begin{gather*}
\delta\vp(h)(u_1,u_2)=\vp(h,u_1)\cdot u_2-\vp(h,u_2)\cdot u_1,\\
Q\vp(h)(u_1,u_2)=\vp(\vp(h,u_1),u_2)-\vp(\vp(h,u_2),u_1)-\vp(h,\th_\m(u_1,u_2)).
 \end{gather*}
For $\nu\in(\La^2\m^*\ot\m)^\h$ let us also define
$p_\nu\in\h^*\ot\La^2\m^*\ot\h$ by the formula
$p_\nu(h)(u_1,u_2)=\vp(h,\nu(u_1,u_2))$ and denote
$\Pi_\vp=\{p_\nu\,\op{mod}B^1(\h,\La^2\m^*\ot\h)\}\subset H^1(\h,\La^2\m^*\ot\h)$.

 \begin{theoremb}\label{newcohomstatement}
The Jacobi identity $\op{Jac}(v_1,v_2,v_3)=0$ with 1 argument from $\h$ and
the others from $\m$ constrains the cohomology $[\vp]\in H^1(\h,\m^*\ot\h)$ so:\\
(1) $[\delta\vp]=0\in H^1(\h,\La^2\m^*\ot\m)$,  whence $\delta\vp=d\theta_\m$;\\
(2) $[Q\vp]\equiv0\in H^1(\h,\La^2\m^*\ot\h)\,\op{mod}\Pi_\vp$, so $Q\vp=d\theta_\h$ for some choices of $\vp,\theta_\m$.
 \end{theoremb}
Now we are ready to proceed with (a sketch of) the proof of Theorem ${\bf3^+}$.
Details of the \textsc{Maple} (DifferentialGeometry) computations can be found in the arXiv supplement.

 \begin{proof}
Let us begin with $\h=\mathfrak{r}$, when this subalgebra contains a grading element $s\in\mathfrak{p}\subset\mathfrak{su}(1,2)$. The cohomology $H^1(\h,\m^*\ot\h)$ has dimension 6, and a 2-dimensional subspace in it satisfies the linear constraint (1) from Theorem \ref{newcohomstatement}. Thus we parametrize $\vp$ by two essential parameters $c_1,c_2$. One may then easily derive the equations $c_1^2=c_2^2=0$ from equation (2) by using that we are working over the real numbers. Hence $[\vp]=0$, and thus $\g$ admits a reductive complement to $\h$, a case which was already completed in \cite{KW1}.

Now let $\h=\mathfrak{l}_1\subset \mathfrak{su}(1,2)$. In this case the cohomology $H^1(\h,\m^*\ot\h)$ has dimension 4. The full cohomology space satisfies equation (1). After solving the linear equation $\delta\vp=d\theta_\m$, the full remaining Jacobi system is simple enough, so that we compute a Groebner basis and obtain all solutions (in this way we get five distinct families of algebras). All the resulting Lie algebras are isomorphic to $\g$ from the example above, and $\mathfrak{l}_1$ has 1-dimensional intersection with $V^4\subset\g$. It is then possible to show that all embeddings of the Lie algebra $\mathfrak{l}_1$ into $\g$ with the required intersection, and which admit an invariant almost complex structure on the quotient $\g/\mathfrak{l}_1$, are conjugate in $\g$. The non-degeneracy of the Nijenhuis tensor can be verified from the structure constants of $\g$ given in the example.

Finally let $\h=\mathfrak{l}_0$. Then the cohomology $H^1(\h,\m^*\ot\h)$ has dimension 10. The solution space of (1) has dimension 5. After solving the linear equation $\delta\vp=d\theta_\m$, the full Jacobi system (and even equation (2) alone) is too complex for Groebner-basis methods or direct solvers, and a full solution seems to be out of reach. However, we compute the Nijenhuis tensor of the corresponding structure $J$ in terms of the Lie algebra parameters, without solving the Jacobi system. The coefficients of this depend only on 10 parameters (out of the 145 variables of the Jacobi system). We then show that modulo the Jacobi system, the Nijenhuis tensor has rank at most 2, hence it is degenerate and we are done.
 \end{proof}



\begin{thebibliography}{WW}
 \footnotesize

\bibitem{AKW}
D.\,V.\ Alekseevsky, B.\,S.\ Kruglikov, H.\ Winther, {\it Homogeneous almost complex structures in dimension 6 with semi-simple isotropy\/}, Ann. Glob. Anal. Geom. {\bf 46}, 361–387 (2014).

\bibitem{KW1}
B. Kruglikov, H. Winther, {\it Almost complex structures in 6D with non-degenerate Nijenhuis tensors
and large symmetry groups\/}, Ann. Glob. Anal. Geom. {\bf 50}, 297–314 (2016).

\bibitem{KW2}
B. Kruglikov, H. Winther, {\it Reconstruction from Representations: Jacobi via Cohomology\/}, arXiv:1611.05334 (2016).

\end{thebibliography}
\end{document}